\documentclass{amsproc}
\usepackage{amssymb,latexsym}

\newtheorem{theorem}{Theorem}[section]
\newtheorem{lemma}[theorem]{Lemma}
\newtheorem{corollary}[theorem]{Corollary}

\theoremstyle{definition}
\newtheorem{definition}[theorem]{Definition}
\newtheorem{example}[theorem]{Example}

\theoremstyle{remark}

\numberwithin{equation}{section}

\newcommand{\Q}{{\mathbb Q}}
\newcommand{\R}{{\mathbb R}}
\newcommand{\N}{{\mathbb N}}
\newcommand{\Z}{{\mathbb Z}}
\newcommand{\T}{{\mathbb T}}
\newcommand{\C}{{\mathbb C}}

\newcommand{\nn}{{\newline}}
\begin{document}

\title{Applications of the work of Stone and von Neumann to wavelets}

\author{Judith A. Packer}
\address{Department of Mathematics, Campus Box 395, University of Colorado, Boulder, CO, 80309-0395}
\email{packer@euclid.colorado.edu}

\date{May 27, 2003}


\begin{abstract}
This survey paper examines the work of J. von Neumann and M.H. Stone as it relates to the abstract theory of wavelets.  In particular, we discuss the direct integral theory of von Neumann and how it can be applied to representations of certain discrete groups to study the existence of normalized tight frames in the setting of Gabor systems and wavelets, via the use of group representations and von Neumann algebras. Then the extension of Stone's theorem due to M. Naimark, W. Ambrose and R. Godement is reviewed, and its relationship to the multiresolution analyses of S. Mallat and Y. Meyer and the generalized multiresolution analyses of L. Baggett, H. Medina, and K. Merrill.  Finally, the paper ends by discussing some recent work due to the author, Baggett, P. Jorgensen and Merrill, and its relationship to operator theory.\end{abstract}

\maketitle


\section{Introduction}

This survey paper is an expanded version of a talk given at the AMS Special Session on Operator Algebras, Quantization, and Noncommutative Geometry: A Centennial Celebration in Honor of J. V. Neumann and M.H. Stone, which was organized at the AMS Annual Meeting in Baltimore by Robert Doran and Richard Kadison, and took place January 16 -- 17, 2003.

In the mid 1980's, what is now known as wavelet theory was developed into a coherent framework, although many of the techniques in this theory had been developed as early as the 1960's. This theory involved the analysis of functions $f$ in $L^2(\mathbb R^n)$ by  studying the frequency content locally in time.  The key theme behind this approach is to analyze a function by studying dilates and translates of an appropriate family of ``wavelets" and/or  ``Gabor functions", which in turn gives either a orthonormal or frame expansion for the function being studied in terms of dilates and translates of a finite number of functions. This theory has led to many important applications in physics, engineering, and signal processing.  In turn, new tools in mathematics, in particular a different approach to Fourier analysis, have been obtained by analyzing the work of physicists and engineering in subdivision schemes and subband coding theory.

It is the point of view of some wavelet theorists that to use the Fourier transform, or more generally, duality for abelian groups, in the study of wavelets, is somewhat circular, because one of the points of wavelet theory is to offer an alternative to ordinary harmonic analysis.  On the other hand, at certain stages in the development of wavelet theory, Fourier analysis has proved critical. This theme is illustrated very clearly by O. Bratteli and P. Jorgensen in their recent book \cite{BJ}. It is the intention of this paper to discuss how certain key aspects of operator theory due to John von Neumann and Marshall Stone arise very naturally in the study of wavelets. Indeed, because wavelet bases and frames and Gabor bases and frames come from dilation, translation, and modulation operators acting on a finite set of functions, it is natural to expect that the spectral analysis of these operators is important in the study of wavelets.  Moreover, dilation, translation, and modulation operators on $L^2(\mathbb R^n)$ are all unitary operators, and various combinations of these operators give rise to countable subgroups of the unitary group of a Hilbert space. 

The approach taken in this article, then, is to think of the key operators in wavelet theory as arising from unitary group representations of countable discrete groups, and then to apply the tools of direct integral theory and von Neumann algebras to their study. Direct integral theory in its simplest form is the analysis of a single unitary operator, and thus reduces to the spectral multiplicity theory of E. Hellinger and H. Hahn, worked out in 1912--1913.  A generalization of this theory led naturally to Stone's theorem concerning the decomposition of a unitary representation of the real line on a Hilbert space by means of a projection-valued measure on $\R$ \cite{S}, which in turn gave a decomposition of the original representation as a direct integral over $\R$ of one-dimensional (hence irreducible) representations. A further generalization of Stone's theorem, due to both Stone and von Neumann, led very naturally to a classification of the irreducible representations of the real three-dimensional Heisenberg group $H_{\R}\;=\;\{\left(\begin{array}{rrr}
1&x&z\\
0&1&y\\
0&0&1
\end{array}\right):\; x,y,z \in\R\}$. 

After a review of direct integral theory for unitary representations, it is our aim to survey how work of Stone and von Neumann, specifically in group representations and direct integral theory, has been used in the study of the theory of wavelets to date. We begin with a discretized version of Stone and von Neumann's approach to the study of representations of the Heisenberg group, that is, with a 
study of unitary representations of the integer lattice of this Heisenberg group,
$H\;=\;\{\left(\begin{array}{rrr}
1&p&r\\
0&1&q\\
0&0&1
\end{array}\right):\; p,q,r \in\Z\}$ that are closely tied up with the theory of Gabor frames; the representation involved generate translation and modulation operators.  We discuss a proof due to L. Baggett \cite{B} of a earlier result on such frames due to I. Daubechies \cite{D1} and M. Rieffel \cite{Ri} which uses direct integral theory. We then move on to a discussion of wavelets and their construction from the wavelet sets of X. Dai and D. Larson \cite{DL}, and discuss a result of L.-H. Lim, J. Packer, and K. Taylor \cite{LPT}.  In this case the group involved will be a discrete semi-direct product group coming from dilation and translation operators on $L^2(\R^n).$
Finally we move on to multi-resolution analyses and some recent work of L. Baggett, H. Medina, and K. Merrill on generalized multiresolution analyses \cite{BMM}.  In these results, the ordinary spectral theory applied to unitary representations of $\Z^n$ plays a key role.  Recent progress in the this area due to Baggett, P. Jorgensen, Merrill and Packer will be outlined \cite{BJMP}.

Naturally in a short paper of this kind much material must be omitted.  In particular, we make no attempt to cover wavelet packets, or any of the exciting work relating group representations of continuous groups to wavelets. We thank Professors Lawrence Baggett, Palle Jorgensen, Wayne Lawton, Kathy Merrill, Arlan Ramsay, Marc Rieffel, Wai-Shing Tang, and Keith Taylor for many helpful conversations on topics related to this paper. 


\section{Direct Integrals, Direct Integral Decomposition and Induced Representations }

 Professor Richard Kadison kindly informed us that the paper of von Neumann ``On rings of operators. Reduction theory" \cite{vN2}, which concerned direct integral theory and its relationship to the decomposition of von Neumann algebras, had been written in the late 1930's but sat in von Neumann's drawer at the Institute of Advanced Study in Princeton until F.J. Mautner asked him about the material in this paper in the late 1940's%
\footnote{Richard Kadison, anecdote given on January 16, 2003, in Baltimore.}%
.  Indeed, von Neumann had published results on the theory of direct integrals of Hilbert spaces as early as 1938 in \cite{vN1}, and he mentions in  the preface \cite{vN2} that the bulk of this ``paper in the drawer" was written during the time period 1937-38. We briefly review the theory in the following paragraphs.  Although our notation is not the same as von Neumann's original notation, the underlying spirit is the same.


\subsection{A review of direct integral theory}

We review the theory of direct integrals of Hilbert spaces and direct integrals of representations. Here we depend heavily on the expositions of J. Dixmier \cite{Dx}, R. Fabec \cite{Fa}, G. Mackey \cite{Mac} and R. Kadison and J. Ringrose \cite{KR}.

Let $(X,{\mathcal B}, \nu)$ be a standard Borel space, that is, ${\mathcal B}$ is the Borel structure of a Polish topology on a topological space $X.$   
\begin{definition}
\label{def hilbbund}
Let ${\mathcal I}$ be a Borel space. We say that ${\mathcal I}$ is a {\bf Hilbert bundle} over $X$ if there is a surjective Borel mapping $p:{\mathcal I}\;\rightarrow\;X$ such that the following conditions are satisfied:
\renewcommand{\labelenumi}{(\roman{enumi})}
\begin{enumerate}
\item For each $x\in X,$ the space $p^{-1}(x)\;=\;{\mathcal H}_x$ has the structure of a Hilbert space.
\item There exists a sequence $\{f_i\}_{i=1}^{\infty}$ of Borel functions $f_i:X\;\rightarrow\;{\mathcal I}$ with $p\circ f_i(x)\;=\; x\;\forall x\in X$ such that the closed linear span of 
$\{f_i(x)\}_{i=1}^{\infty}$ in ${\mathcal H}_x$ is equal to ${\mathcal H}_x$ for all $x\in X.$
\item For all $i,j\;\in \N,$ the function $x\mapsto <f_i(x),f_j(x)>_{{\mathcal H}_x}$ is Borel.
\item The function $g:X\;\rightarrow\;{\mathcal I}$ is Borel if and only if $p\circ g$ is Borel and the functions 
$$x\mapsto\;<g(x),f_i(p(g(x)))>{{\mathcal H}_{p(g(x))}}$$
are Borel for all $i.$ 
\end{enumerate}
\end{definition}

The {\bf Borel cross-sections} to the Hilbert bundle are those Borel functions $f:X\;\rightarrow\;{\mathcal I}$ which satisfy $p\circ f(x)\;=\;x.$ 
\begin{definition}
\label{def dirinth}
A Borel cross-section $f:X\;\rightarrow\;{\mathcal I}$ is said to be square-integrable if $\int_X||f(x)||^2d\nu <\infty.$  The set of square-integrable  cross sections is an inner product space, with the inner product $\langle f, g\rangle$ given by the obvious formula
$$ \langle f, g\rangle\;=\;\int_X <f(x), g(x)> d\nu.$$
 Upon identifying  square-integrable  cross sections $f,f'$ satisfying $\int_X||f(x)-f'(x)||^2 d\nu = 0,$ we obtain the {\it direct integral Hilbert space} $\int_X^{\oplus}[{\mathcal H}_x]d\nu$.
\end{definition} 

\begin{example}
Let ${\mathcal H}$ be a fixed separable Hilbert space, and suppose that ${\mathcal I}\;=\;X\times {\mathcal H}.$ Define $p:{\mathcal I}\rightarrow X$ to be projection in the first variable. Then ${\mathcal H}_x\;=\;{\mathcal H},\;\forall x\in {\mathcal H}.$ Fixing a basis $\{e_i\}\;\subseteq\;{\mathcal H}$ and defining $f_i(x)=e_i\forall x\in X,$ we see that ${\mathcal I}$ is a Hilbert bundle.
 The vector space of all Borel functions 
$f:X\;\rightarrow\;{\mathcal H}$ such that $<f(x),e_i>$ is measurable for every $e_i$ gives the Borel cross-sections to the bundle. When we consider the direct integral Hilbert space $\int_X^{\oplus}[{\mathcal H}_x]d\nu$ in this case, we obtain exactly the family of square-integrable functions from $X$ into ${\mathcal H}$, where two functions are identified if they are equal except on a set of $\nu$-measure $0$.  The associated Hilbert bundle is called a {\bf constant bundle}, and the direct integral is frequently denoted by $L^2(X,{\mathcal H})\;=\;L^2(X)\otimes {\mathcal H}.$  Specializing further to the case ${\mathcal H}=\C,$ we obtain exactly $L^2(X).$ 
\end{example}

We note that by condition (ii) of Definition \ref{def hilbbund}, each Hilbert space ${\mathcal H}_x$ is separable.  Define a functions $\text{dim}:X\;\rightarrow\; \N\cup\;\{\infty\}$ by $\text{dim}(x)\;=\;\text{dim}({\mathcal H}_x).$ For each $j\in\N\cup\{\infty\}$ let 
$X_j\;=\;\text{dim}^{-1}(j).$ 
Von Neumann proved that for any Hilbert bundle defined as above, and for each $j\;\in\N\cup\{\infty\},\;X_j$ is measurable. 
Hence the space $X$ can be partitioned as the disjoint union of measurable subsets 
$$X\;=\;\bigsqcup_{j\in\N\cup\{\infty\}}X_j,$$
and it is known that the direct integral Hilbert space $\int^{\oplus}_{X_j}[{\mathcal H}_x]d\nu_{|X_j}$ is isomorphic to the direct integral of the constant bundle ${\mathcal H}(j),$ where 
${\mathcal H}(j)$ represents a fixed separable Hilbert space of dimension $j$. 

Given a direct integral Hilbert space $\int_X^{\oplus}[{\mathcal H}_x]d\nu,$ we recall how one forms direct integrals of bounded operators to act on the Hilbert space.
Suppose we are given a Hilbert bundle ${\mathcal I},$ and suppose that $x\mapsto T(x)$ is a function that assigns to each $x\in X$ an element of $B({\mathcal H}_x),$ the algebra of all bounded operators on ${\mathcal H}_x.$ 
\begin{definition}
\label{def dirintop}
Let $({\mathcal I}, p, X)$ be a Hilbert bundle.
We say that the mapping $x\mapsto T(x)$ is a {\it measurable cross-section of bounded operators} if for every Borel cross-section $f:X\;\rightarrow {\mathcal I},$ the cross section $x\mapsto\;T(x)f(x)$ is also in ${\mathcal I}.$
If $x\mapsto T(x)$ is a measurable cross-section of operators, we say it is a {\it bounded} measurable cross-section if the function $x\mapsto \|T(x)\|$ is in $L^{\infty}(X).$ In this case, 
$x\mapsto T(x)$ determines an element of ${\mathcal B}(\int_X^{\oplus}[{\mathcal H}_x]d\nu),$ denoted by
$$T\;=\;\int^{\oplus}_X\;[T(x)]d\nu,$$ where we have 
$T(f)(x)\;=\;T(x)(f(x))$ for any square-integrable  cross sections $f\in\;\int_X^{\oplus}[{\mathcal H}_x]d\nu.$
Conversely, we say a bounded operator $T\in B(\int_X^{\oplus}[{\mathcal H}_x]d\nu)$ is {\it decomposable} if there is a bounded measurable cross-section $x\mapsto T(x)$ such that 
$$T\;=\;\int^{\oplus}_X\;[T(x)]d\nu.$$
\end{definition}

Now we briefly discuss the direct integral decomposition of unitary group representions.  Suppose that $G$ a locally compact second countable group.  Let $U:G\;\rightarrow\;{\mathcal U}({\mathcal H})$ be a continuous unitary group representation of $G,$ that is, a homomorphism of the group $G$ into the unitary group ${\mathcal U}({\mathcal H}),$ such that $\forall x,y\in {\mathcal H},$
the map 
$$g\;\mapsto\;<U_g(x),y>$$
is continuous from $G$ to $\C.$ 
\begin{definition}
\label{def dirintrep}
We say that $U$ is decomposable as a direct integral over the Borel space $(X,\nu)$ if the following conditions hold:
\renewcommand{\labelenumi}{(\roman{enumi})}
\begin{enumerate}
\item The Hilbert space ${\mathcal H}$ is isomorphic to a direct integral $\int_X^{\oplus}[{\mathcal H}_x]d\nu.$
\item For each $x\in X$ there is a unitary group representation $U(x):G\;\rightarrow\;{\mathcal H}_x$ such that 
$$U_g\;=\;\int_X^{\oplus}[U(x)_g]d\nu,\;\forall g\in G,$$
with respect to the identification of ${\mathcal H}$ with $\int_X^{\oplus}[{\mathcal H}_x]d\nu$
made in (i).
\end{enumerate}
\end{definition}

The theory of von Neumann algebras and in particular, the paper \cite{vN2}, was fundamental in showing that every unitary decomposition of $G$ could be decomposed as a direct integral of irreducible representations of $G.$ 
In order to briefly describe this procedure, we first recall the definition of a von Neumann algebra ${\mathcal M}$ of 
${\mathcal B}({\mathcal H}),$ due to F. Murray and J. von Neumann \cite{MvN1,MvN2}.  As other articles in this survey article will go into these definitions in more detail, here we deliberately make our remarks sketchy and refer the reader to \cite{Dx}, for example, for full details.
\begin{definition}
\label{def vNalg}
A {\bf von Neumann algebra} on the Hilbert space ${\mathcal H}$ is a subalgebra ${\mathcal M}$ of ${\mathcal B}({\mathcal H})$ that is closed under the adjoint operation, contains the identity operator $\text{Id}_{{\mathcal H}},$ and is closed in the strong operator topology. The set of all operators in ${\mathcal B}({\mathcal H})$ that commute with every element of ${\mathcal M}$ is called the {\bf commutant} of ${\mathcal M}$ and is denoted by ${\mathcal M}'.$ 
\end{definition}
Murray and von Neumann proved that if ${\mathcal M}$ is a von Neumann algebra, then its commutant ${\mathcal M}'$ is a von Neumann algebra also, and 
$$[{\mathcal M}']'\;=\;{\mathcal M}.$$
Indeed, if ${\mathcal N}$ is any subalgebra of ${\mathcal B}({\mathcal H})$ which is closed under the adjoint operation, they proved that the smallest von Neumann algebra in ${\mathcal B}({\mathcal H})$ containing ${\mathcal N}$ is 
$[{\mathcal N}']'.$ For this reason the double commutant ${\mathcal N}''$ is called the von Neumann algebra {\bf generated by ${\mathcal N}.$}
A von Neumann algebra ${\mathcal M}$ is said be a {\bf factor} if 
$${\mathcal M}\;\cap\;{\mathcal M}'\;=\;\{\lambda \text{Id}:\;\lambda\in\C\}.$$
This means exactly that the center of ${\mathcal M}$ consists of scalar multiples of the identity operator.  Von Neumann and Murray showed that one can define a dimension function, which is unique up to a normalizing scalar, on projections of a factor ${\mathcal M},$ which extends to a trace function on the factor ${\mathcal M}.$  A factor ${\mathcal M}$ is said to be 
of type $II_1$ if  this dimension function takes on exactly the values $[0,1].$

Of course not all von Neumann algebras are factors, and there is a vast literature describing the structures and decomposition of von Neumann algebras.  Here we only touch on that part of von Neumann algebra theory relevant to our discussion.  We first use von Neumann algebras to describe how to 
decompose unitary representations of groups into direct integrals.  This procedure is due to von Neumann \cite{vN2}, F. Mautner \cite{Mau} and R. Godement \cite{Go}, and goes as follows:
\renewcommand{\labelenumi}{(\roman{enumi})}
\begin{enumerate}
\item Form ${\mathcal G}_U\;=\;\{U_g:g\in G\}'',$ the von Neumann subalgebra of $B({\mathcal H})$ generated by $\{U_g:g\in G\};$
\item Let ${\mathcal M}$ be a maximal abelian subalgebra of the commutant of ${\mathcal G}_U,$ and find an isomorphism between ${\mathcal M}$ and $L^{\infty}(X,\nu)$ for some measure space 
$(X,\nu);$
\item Decompose ${\mathcal H}$ as $\int_X^{\oplus}[{\mathcal H}_x]d\nu$ and $U$ as $\int_X^{\oplus}[U(x)_g]d\nu.$ 
\end{enumerate}
 However, since the maximal abelian subalgebras in ${{\mathcal G}_U}'$ are not unique up to inner conjugacy, this direct integral decomposition need of $U$ not be unique, and in the case where $G$ is countable and discrete and $[G: [G,G]]$ is infinite, there always exist representations having non-uniqueness in their direct integral decompositions; see \cite{Mac}, Chapter 3, Section 5 for an example.


\subsection{A review of induced representations}

Induced representations will play an important role in the discussion that follows, 
so we review their definition at this point.  We follow the description of induced representations given
in Chapter X of \cite{K}. 
\begin{definition}
Let $G$ be a locally compact second countable 
group with closed subgroup $N$, and suppose we are given a Borel cross section
$s:N\backslash G \to G$ with $s([e_G])=e_G,$ where here $e_G$ represents the identity element of the group $G.$ 
Use $s$ to define a one-cocycle $\sigma$ for the
right action of $G$ on the coset space $N\backslash G$ as follows:  
$\sigma: N\backslash G\times G \to N$ is defined by
$$\sigma(x, g)\;=\;s(x)g [s(x \cdot g)]^{-1},\;(x,g)\in\;N\backslash G\times G.$$ 
Let $\rho$ be a unitary representation 
of the group $N$ on the Hilbert space $\mathcal{H}$. Then the induced
representation $Ind_{G_0}^G(\rho)$ has as its representation space
$L^2(N\backslash G,
\mathcal{H}).$
where the measure $\nu$ on $N\backslash G$ is quasi-invariant
under translation by $G$. If this measure is invariant, then the formula
for $Ind_{N}^G(\rho)$ is given by
$$Ind_{N}^G(\rho)(g)f](x) =\rho(\sigma(x,g))[f(xg)],$$
for $g\in G,\;x\in N\backslash G,\;f\in L^2(N\backslash G,
\mathcal{H})$. 
\end{definition}
In the discussion that follows, all of the groups we study will be discrete, and the
translation-invariant measure on the coset space $N\backslash G$
will just be the counting measure.

\begin{example}
\label{ex heis}
Let $H$ denote the integral Heisenberg group; that is, $H$ consists of all matrices of the form
$$\{\left(\begin{array}{rrr}
1\;&p\;&\;r\\
0\;&1\;&q\;\\
0\;&0\;&1
\end{array}\right)\;:\;p,q,\;r\;\in\;\mathbb Z\}.$$
It is easier notationally to parametrize $H$ by the set $\{(p,q,r)\;:\;p,q,r \in \Z\},$ with obvious operations coming from the matrix multiplications. Let $N$ be the abelian subgroup of $H$ defined by $N\;=\;\{(0,q,r)\; ;\;q,r \in \Z\}.$ Fix $t,u\;\in\;[0,1)$ and define a one-dimensional representation or character $\chi(u,t)$ of $N$ on $\C$ by 
$$\chi(u,t)_{(0,q,r)}\;=\;e^{2\pi i qu}e^{2\pi i rt},\;(0,q,r)\;\in\;N.$$
The coset space $H_0\backslash H$ can be identified with $\Z,$ and the cross-section 
$s:\; N\backslash H\;\cong\;\Z\;\rightarrow\;H$ is given by 
$$s(k)\;=\;(0,0,k),\;k\in\;N\backslash H.$$
The corresponding one-cocycle for the action of $H$ on $N\backslash H$ is given by
$$\sigma(k,(p,q,r))\;=\;(p,q,0).$$
Hence $Ind_{N}^H(\chi(u,t))$ is defined on $l^2(\Z,\C)\;\cong\;l^2(\Z)$ by 
$$Ind_{N}^H(\chi(u,t))(p,q,r)f(k)\;=\;e^{2\pi i qu}e^{2\pi i rt}f(k+r),\;f\in l^2(\Z),\;(p,q,r)\in\; H.$$  
We call the representations $Ind_{N}^H(\chi(u,t))$  {\bf ``Stone-von Neumann"} representations of $H,$ because they can be extended to the three-dimensional real Heisenberg group and for $t\not= 0$ provide irreducible representations of the real Heisenberg group of the type first discovered by Stone and von Neumann. We remark that elsewhere in the literature these representations are referred to as {\bf ``Weyl-Heisenberg"} representations, since W. Heisenberg first described commutation relations between the position and momentum operators, and H. Weyl, using the work of von Neumann and Stone, noticed that by exponentiating these operators,  one could construct unitary operators on a Hilbert space satisfying a key commutation relation that will arise later in this paper.
\end{example}

We note that the operations of taking direct integrals of representations and inducing representations commute:  That is, if $N$ is a closed subgroup of the locally compact second countable group $G,$ and if $L$ is a representation of $N$ which is written as a direct integral 
$$L_n\;=\;\int^{\oplus}_X [L_n(x)]d\nu,\;n\in N$$
on the Hilbert space $\int^{\oplus}_X[{\mathcal H}_x]d\nu,$
then the induced representation $Ind_{N}^G(L)$ is unitarily equivalent to the direct integral of representations 
$$\int_X^{\oplus}[Ind_{N}^G(L_x)]d\nu$$ defined on the direct integral Hilbert space 
$$\int_X^{\oplus}[L^2(N\backslash G,{\mathcal H}_x)]d\nu.$$
We refer the reader to \cite{Mac} and \cite{Fa} for details.

 \section{Important operators and groups for wavelet and frame theory and the von Neumann algebras they generate}

We first review the notion of a normalized tight frame in a Hilbert space.
\begin{definition} 
  A sequence $\{ f_n : n \in {\mathbb N} \}$ of elements in a Hilbert
  space $\mathcal H$ is said to be a {\it frame} if there are real
  constants $C,D > 0$ such that
  \begin{equation} 
\label{eq-ntframe}
  C\sum_{n=1}^{\infty} |\langle f,f_n \rangle|^2\;\leq\; \|f\|^2 \;\leq\;
      D\sum_{n=1}^{\infty} |\langle f,f_n \rangle|^2   
\end{equation}
  for every $f \in \mathcal H$.  The constants $C$ and $D$ are called the {\bf frame bounds} for the frame. If $C=D,$ the frame is called a {\bf tight frame}, and if $C=D=1,$ the frame is called a {\bf normalized tight frame}, abbreviated NTF. It is known that $\{ f_n : n \in {\mathbb N} \}$ is an NTF for ${\mathcal H}$ if and only if the following reconstruction formula holds: 
\begin{equation} \label{recon}
f\;=\;\sum_{n=1}^{\infty}\langle f,f_n \rangle f_n,\;\forall f\in {\mathcal H}.
\end{equation} 
\end{definition}
Because of this last identity, called the {\bf reconstruction formula}, normalized tight frames are sometimes referred to as ``Parseval frames" in the literature. Clearly, an orthonormal basis for ${\mathcal H}$ is a special case of a normalized tight frame.


\subsection{Key operators in Gabor and wavelet frame theory}

It turns out that in the construction of Gabor and wavelet frames, there are three key types of unitary operators that play important roles. They are dilation operators, translation operators, and modulation operators, which are defined on $L^2(\R^n).$ 
\begin{definition}
\label{def dil}
 Let $A$ be an $n\times n$ integral matrix all of whose eigenvalues have modulus greater than one.
The dilation operator corresponding to $A$ is defined on $L^2(\R^n)$ by 
$$D_A(f)(t)\;=\;|\text{det}(A)|^{\frac{1}{2}}f(At).$$
\end{definition}
\begin{definition}
\label{def trans}
Let $\Gamma$ be a cocompact lattice in $\R^n,$ for example, $\Gamma = \Z^n,$ and let $\gamma\in \Gamma.$
The translation operator corresponding to $\gamma$ is defined on $L^2(\R^n)$ by 
$$T_{\gamma}(f)(t)\;=\;f(t-\gamma).$$
\end{definition}
\begin{definition}
\label{def mod}
Let $\alpha\;\in\;\R^n.$  The modulation operator corresponding to $\alpha$ is defined by 
$$M_{\alpha}(f)(t)\;=\;e^{2\pi i\alpha\cdot t}f(t).$$
\end{definition}
Let ${\mathcal F}$ denote the Fourier transform on $L^2(\mathbb R),$ defined by 
$${\mathcal F}(f)(x)\;=\;\int_{\mathbb R}f(t)e^{-2\pi ix\cdot t}dt.$$  
Then one can verify the following relations:
$$T_{\gamma}\;D_A \;=\;D_AT_{A(\gamma)},$$
$$M_{\alpha}T_{\beta}\;=\;e^{-2\pi i \alpha \beta}T_{\beta}M_{\alpha},$$ 
$${\mathcal F}D_A{\mathcal F}^{-1}\;=\;D_{{A^T}^{-1}},$$
and 
$${\mathcal F}T_{\gamma}{\mathcal F}^{-1}\;=\;M_{\gamma}.$$

We now use the operators listed above to define certain kinds of frames in $L^2(\R^n).$
\begin{definition} 
\label{def Gaborframe}
Fix full-rank lattices ${\mathcal L}$ and ${\mathcal K}$ in $\R^n.$ By this we mean that ${\mathcal L} = L(\Z^n)$ and ${\mathcal K}\;=\;K(\Z^n)$ for two non-singular $n\times n$ real matrices $L$ and $K.$  
We say that a vector $f\;\in L^2(\mathbb R^n)$ forms a Gabor type NTF for $L^2(\R^n)$ 
(respectively, a Gabor orthonormal basis for $L^2(\R^n)$) with respect to the lattices ${\mathcal L}$ and ${\mathcal K}$ if the set of functions 
$$\{M_{\lambda}T_{\kappa}(f)\;:\;\lambda\in {\mathcal L},\;\kappa\in {\mathcal K}\}$$
forms an NTF for  for $L^2(\R^n)$ (respectively, a orthonormal basis for $L^2(\R^n)$).
\end{definition}
Specializing to the case $n=1,$ we say that a vector $f\;\in L^2(\mathbb R)$ forms a Gabor NTF for $L^2(\R)$ with respect to the constants $\alpha,\;\beta>0$ if 
the set of functions 
$$\{M_{p\alpha}T_{\beta q}(f)\;:\;p,q\;\in\;Z^n,\;1\;\leq i\; m\}$$
forms an NTF for  for $L^2(\R).$ 

\begin{example}
\label{ex Gaborchac}
Let $n=1$ and $\alpha\;\beta\;=\;1,$ and take $f(t)=\chi_{[0,1]}(t).$  An easy calculation with Fourier series shows that $\{M_q T_p(f)\;:\;p,q \in \Z\}$ is an orthonormal basis for $L^2(\mathbb R),$ hence is a Gabor type NTF.
\end{example}
\begin{example} 
 Again take $n=1$ and $\alpha\;\beta\;=\;1.$ Let $f$ be the Gaussian $f(t)=\;e^{-t^2},$ then the set of functions $\{M_q T_p(f)\;:\;p,q \in \Z\}$ form a set whose finite linear combinations are dense in  $L^2(\mathbb R).$ 
Von Neumann knew this was true in the 1940's. However, in this case $\{M_q T_p(f)\;:\;p,q \in \Z\}$ is not a frame, c.f. \cite{Ben},\cite{D1}. 
\end{example} 

\begin{definition}
\label{def wavelet}
A {\bf {NTF wavelet family }} for dilation by $A$ is a subset $\{\psi_1,\cdots,\psi_m\}\subseteq L^2(\mathbb R^n)$ such that
$$\{D_A^j T_v(\psi_i):\;1\leq i \leq m,\;j\in\mathbb Z,\;v\in\mathbb Z^n\}$$
is an NTF for $L^2(\mathbb R^n).$ If the functions $\{D_A^j T_v(\psi_i):\;1\leq i \leq m,\;j\in\mathbb Z,\;v\in\mathbb Z^n\}$ form an orthonormal basis for $L^2(\R^n),$ the family $\{\psi_1,\cdots,\psi_m\}$ is called an orthonormal wavelet family for dilation by $A.$
\end{definition}
Thus when one has a NTF wavelet family $\{\psi_1,\cdots,\psi_m\}$ for dilation by $A,$ for every $f\in\;L^2(\R^n)$ the reconstruction formula Equation \ref{recon} is satisfied with respect to the NTF $\{D_A^j T_v(\psi_i):\;1\leq i \leq m,\;j\in\mathbb Z,\;v\in\mathbb Z^n\}.$

In applications, the function $f$ being studied in $L^2(\mathbb R^n)$ determines the wavelet family $\{\psi_1,\cdots,\psi_m\}$ selected, and one would hope that all but finitely many of the wavelet coefficients in the reconstruction formula for $f$ would be zero or at least so small as to be safely ignored in computations, without loss of essential data. 
\begin{example}
\label{classical case}
Let $n=1$ and $A=2,$ and define
$$\psi_H(t)=\left\{\begin{array}{ll}
1&\mbox{if}\ \;x\in [0,\frac{1}{2}],\\
-1 &\mbox{if}\ \;x\in (\frac{1}{2}, 1]\\
0&\;\mbox{otherwise}.
\end{array}\right.$$
Then $\{D_2^j T_v(\psi_H):\;j,v\;\in\Z\}$ forms an orthonormal basis for $L^2(\R),$ and $\psi_H$ is called the {\bf Haar wavelet} for dilation by $2.$
\end{example}

\subsection{Group representations associated to Gabor and wavelet theory}

The operators connected to Gabor frames and wavelet frames are very much connected to 
certain groups, as we describe now.
\begin{definition} 
Let $H$ denote the integral Heisenberg group defined in Example \ref{ex heis}.
Define $U^{\alpha,\beta}_{(p,q,r)}(f)(t)=e^{-2\pi i r\alpha \beta}e^{2\pi iq\alpha t}\;f(t-p\beta).$ 
Then  $U^{\alpha,\beta}$ is a unitary representation of $H$ on $L^2(\mathbb R).$ A simple change of variable in $t$ shows that $U^{\alpha,\beta}$ is unitarily equivalent to $U^{\alpha\beta,1}$ so depends only on the product $\alpha\beta$ and not on the separate parameters $\alpha$ and $\beta.$ We write $U^{\delta}$ for $U^{\delta,1}.$
\end{definition}
Note that 
$$U^{\delta}(p,0,0)\;=\;T_{p},$$  
and 
$$U^{\delta}(0,q,0)\;=\;M_{q\delta}.$$
It follows that $f\;\in\;L^2(\R)$ gives an Gabor NTF (respectively, a Gabor orthonormal basis) for 
$L^2(\R)$ if and only if $\{U^{\delta}(p,q,0)(f): p,q\in\Z\}$ forms an NTF (respectively, an orthonormal basis) for $L^2(\R).$  This description can be extended full-rank lattices in the $2n+1$ dimensional real Heisenberg group as well, and we leave the calculations to the reader.
Thus the study of Gabor frames is intimately connected to the study of the representation $U^{\delta}$ of $H$ and to the von Neumann algebra generated by the operators 
$\{U^{\delta}(p,q,r):\;(p,q,r)\in H\}.$ In turn, these representations of $H$ can be viewed as 
projective (or ``multiplier" representations) of $\Z^2$ corresponding to the multiplier $\omega_{\delta}$ on $\Z^2$ defined by 
$$\omega_{\delta}((p_1,q_1),(p_2,q_2))\;=\;e^{-2\pi i p_1q_2\delta},\;(p_1,q_1),(p_2,q_2)\;\in\;\Z^2.$$
Such representations have been studied extensively in the literature, c.f. \cite{bagklep,Rie}. 

We now discuss discrete groups essential to the study of wavelets. We note that for the classical case where $n=1$ and $A=2,$ it was X. Dai and D. Larson and, independently,  B. Brenken, who first realized the connection between the study of representations of this group and wavelet theory, c.f. \cite{DL,Br}.
\begin{definition}
Let $A$ be an integral $n\times n$ dilation matrix, and let 
$\Q_A$ denote the discrete $A$-adic rational subgroup of $\Q^n$ defined by 
$$\Q_A\;=\;\cup_{j\in\Z}[A^{-j}(\Z^n)].$$
The matrix $A$ defines an action of $\Z$ on ${\Q}_A$ by automorphisms which we denote by $\theta:$
$$\theta(k)(\beta)\;=\;A^{-k}(\beta),\;\beta\in\;\Q_A,\;k\in\Z.$$
We thus can form the semidirect product of $\Q_A$ and $\Z,\;
\Q_A\rtimes_{\theta}\Z,$ which we call the {\it wavelet group} (for dilation by $A$) and denote by $G_A.$
\end{definition}
We parametrize $G_A$ by pairs $\{(\beta,k)\;:\;\beta\in \Q_A,\; k\in \Z\},$ where the group operation
is defined by 
$$(\beta_1,k_1)\cdot (\beta_2,k_2)\;=\;(\beta_1 +\theta(k_1)(\beta_2), k_1+k_2).$$
$G_A$ is sometimes called a {\it generalized Baumslag-Solitar group}, c.f. \cite{MV}.
\begin{definition}
\label{def waverep}
Define a group representation $W$ of the wavelet group $G_A$ on $L^2(\R^n)$ by
$$W_{(\beta,k)}(f)(t)\;=\;|\text{det}(A)|^{\frac{k}{2}}f(A^k(t-\beta)),\;(\beta,m)\in G_A.$$
We call this the wavelet representation of $G_A.$
\end{definition}
One calculates that 
$$W_{(\beta,k)}\;=\;T_{\beta}D_A^k.$$
From this observation we see that $\{\psi_1,\cdots, \psi_m\}$ is a NTF wavelet family (respectively,  orthonormal wavelet family) for dilation by $A$ if and only if 
$$\{W_{(\theta(k)(\beta),k)}(\psi_i)\;:\beta\in \Z^n,\; k\in\Z,\; 1\leq i\leq m\}$$ forms a NTF
(respectively, an orthonormal basis) for $L^2(\R^n).$
Again it follows that the study of wavelets is intimately connected to the study of the representation $W$ and the von Neumann algebra generated b the unitary operators $\{W_{(\beta,k)}\;:\;(\beta,k)\in G_A\}$. 

 Much of the remainder of the paper will be devoted to the decomposition of the representations $U^{\delta}$ and $W$ via direct integral theory.


\section{The use of direct integrals in the study of Gabor frames and wavelets}

We now use the representations described above in the study of properties of group representations and wavelets.


\subsection{The study of Gabor frames using direct integral theory and von Neumann algebras}

Fix $\delta \;>\;0,$ and consider the Gabor representation $U^{\delta}$ defined in the previous section.
The following decomposition theorem is explicitly due to L. Baggett, although various people had been aware of it previously; we discuss higher dimensions at the end of this section.
\begin{theorem}
\label{thm bag}
\cite{B}
Let $H$ and $N$ denote the integer Heisenberg group and its doubly generated abelian subgroup 
$N\;=\;\{(0,q,r)\; ;\;q,r \in \Z\}$ defined in \ref{ex heis}. Let $E$ denote the interval $[0,\delta)$ and let 
$\nu$ denote Haar measure on $E.$ 
\renewcommand{\labelenumi}{(\roman{enumi})}
\begin{enumerate}
\item There is a unitary equivalence 
$${\mathcal U}:L^2(\R)\;\rightarrow\; \int_E^{\oplus}[l^2(\Z)] d\nu$$
and a one parameter family of representations $\{U^{\delta}(t): t\in E\}$ of $H$ on $l^2(\Z)$ such that $U^{\delta}$ is unitarily equivalent via ${\mathcal U}$ to the direct integral representation 
$$\int_{E}^{\oplus}[U^{\delta}(t)] d\nu$$ of $H$ on $\int_E^{\oplus}[l^2(\Z)] d\nu.$
For $t\in E,$ the representation $U^{\delta}(t)$ is defined by 
$$U^{\delta}(t)_{(p,q,r)}g(k)\;=\; e^{-2\pi i\delta r}e^{-2\pi i qt} e^{2\pi i q\delta k}g(k-p),\;g\in l^2(\Z).$$
\item For each $t\in E,$ let $\chi_{t,\delta}$ denote the one-dimensional representation of $N$ on $\C$ defined in \ref{ex heis}, i.e. 
$$\chi_{t,\delta}((0,q,r))\;=\;e^{-2\pi i\delta r}e^{-2\pi i qt},\;(0,q,r)\;\in\;N.$$
Then each representation $U^{\delta}(t)$ of $H$ is unitarily equivalent to the monomial representation 
$\text{Ind}_{N}^H(\chi_{t,\delta})$ of $H$ on $l^2(N\backslash H)$ induced from the one-dimensional representation $\chi_{t,\delta}$ of $N$ on $\C.$  
\end{enumerate}
\end{theorem}
It follows from the above theorem that the study of the monomial Stone-von Neumann representations $U^{\delta}(t)$ of $H$ are key to understanding Gabor frames. 
It is possible to use this direct integral decomposition to approach the following question. 
\nn {\bf Question:} Does there exist  a single function $f\;\in\;L^2(\mathbb R)$ such that the family
$$\{f_{p,q}(t)\;=\;e^{2\pi iq\alpha t}\;f(t-p\beta):\;p,q\;\in\;\mathbb Z\},$$ or equivalently, the family 
$$\{U^{\delta}_{(p,q,0)}(f):\; p,q\;\in\;\Z\},$$
forms a NTF for $L^2(\mathbb R)$?
\nn {\bf Answer of I. Daubechies \cite{D1}:}  (1) If $\delta\;=\;1$, the answer to the question is yes (Example \ref{ex Gaborchac} gives an orthonormal basis $\{f_{p,q}\}$ ); (2) if $\delta\;<\;1,$ the answer is yes, (3) if $\delta\;>\;1,$ the answer is {\bf no}.

The first proof of (3) for $\delta>1$ and rational is due to Daubechies \cite{D1}. The proof of (3) for $\delta>1$ and irrational is due to Rieffel who used von Neumann algebra ``coupling constant" techniques described in \cite{Ri}. He related this work to his study of the irrational rotation algebra in \cite{Rie}, as certain representations of $H$ give rise to this well-known $C^{\ast}$-algebra. We sketch Rieffel's proof in later paragraphs, but we first describe how this problem has been studied by using direct integral theory. In 1990 L. Baggett gave a proof of the same results involving the explicit direct integral decomposition of the Stone-von Neumann representation mentioned above. 

\begin{example}  
\label{ex Gabor}
We briefly discuss the above result from the direct integral point of view:
for simplicity, we consider $L^2(\R)$ within the framework of Baggett's direct integral decomposition 
$$L^2(\R)\;\cong\;\int_{[0,\delta)}^{\oplus}[l^2(\Z)] d\nu,$$ and consider for what 
values of $\delta$ the function $$f_{\delta}=\chi_{[0,\delta)}\equiv\;\int_{[0,\delta)}^{\oplus}[\delta_{0}]d\nu$$ can
generate a tight frame under the orbit of the operators $$\{\int_{[0,\delta)}^{\oplus}U^{\delta}(t)_{(p,q,0)}d\nu\}.$$ 
(Here $\delta_0$ represents the Dirac delta function at $0$ defined on $\Z.$)

We note that 
$$\{U^{\delta}_{(p,q,0)}(f_{\delta}):p,q\in\Z\}
\;=\;\{\int_{[0,\delta)}^{\oplus}U^{\delta}(t)_{(p,q,0)}f_{\delta}(k)d\nu: p,q, k\in\Z\}$$
$$=\; \{\int_{[0,\delta)}^{\oplus}e^{-2\pi i qt} e^{2\pi i q\delta k}f_{\delta}(k-p): p,q, k\in\Z\}\;=\;\{\int_{[0,\delta)}^{\oplus}e^{-2\pi i qt} e^{2\pi i q\delta k}\delta_{p}(k): p,q, k\in\Z\}.$$
When we consider the direct integral as a tensor product $L^2[0,\delta)\otimes l^2(\Z),$
this last set becomes 
$$\{e^{2\pi i q\delta k}e^{-2\pi i qt}\otimes\delta_{p} : p,q\in\Z\}.$$ 
Let $\delta\leq\;1.$ From our above discussion, one sees that    
$\{U^{\delta}_{(p,q,0)}(f_{\delta}):\; p,q\in\Z\}$ will form a frame for
$L^2(\R)$ if and only if the family of functions $\{e^{2\pi i qt}: \;
q\in \Z \}$ forms a frame for $L^2[0,\delta).$ Let us consider the normalized family 
$\{\frac{1}{\sqrt{\delta}}e^{2\pi i qt}: \;
q\in \Z \}$ for $\delta\leq 1.$
We take $f\in L^2[0,\delta)$ and compute 
$$\sum_{q\in\Z}|\langle f, \frac{1}{\sqrt{\delta}}e^{2\pi i qt}\rangle_{L^2[0,\delta)}|^2$$
$$\;=\;\sum_{q\in\Z}|\int_{[0,\delta)}f(t)\overline{\frac{1}{\sqrt{\delta}}e^{2\pi i qt}}d\nu(t)|^2$$
$$\;=\;\frac{1}{\delta}\sum_{q\in\Z}|\int_{[0,1)}f(t)e^{-2\pi i qt}d\nu(t)|^2$$
(since $f$ is supported on $[0,\delta)$)
$$\;=\;\frac{1}{\delta}\|f\|_{L^2[0,1)}^2$$
(by Parseval's formula)
$$\;=\;\frac{1}{\delta}\|f\|_{L^2[0,\delta)}^2$$
(since $L^2[0,\delta)\subseteq L^2[0,1)$).
Thus if $0<\delta\;\leq 1$ we obtain a tight frame for $L^2[0,\delta)\otimes l^2(\Z)$ with frame bound $\delta.$ Indeed if $\delta\;=\;1$ we get an orthonormal basis via this method, as we have already seen in Example \ref{ex Gaborchac}.

On the other hand, if $\delta>1$ the set $\{e^{2\pi i qt}:\;q\in\Z\}$ does not even span $L^2[0,\delta),$ so one cannot hope to get a frame from $\{U^{\delta}_{(p,q,0)}(f_{\delta}):p,q\in\Z\}$ in this case.
\end{example}

We now discuss this result from the point of view of M. Rieffel, who gave the first proof that if $\alpha\cdot\beta$ was irrational and greater than $1,$ it was impossible for a vector $f\;\in L^2(\mathbb R)$ to form a Gabor NTF in $L^2(\R)$ with respect to the constants $\alpha,\;\beta.$ \cite{Ri}. As before, set $\delta =\alpha\cdot\beta.$ A special case of a result of M. Takesaki \cite{Ta} shows that is possible to form the finite transformation group von Neumann algebras $M(\delta\Z, \R/\Z)$ acting on $L^2(\R)$ on the left and $M(\Z,(\delta\Z)\backslash \R)$ acting on $L^2(\R)$ on the right in such a way that each von Neumann algebra is the other's commutant. [We note that it is now standard to write such transformation group von Neumann algebras as $L^{\infty}(\R/\Z)\rtimes \delta \Z$ and $L^{\infty}((\delta\Z)\backslash \R)\rtimes\Z$ respectively].  Rieffel showed that $M(\Z,(\delta\Z)\backslash \R)$ is generated by two operators $U$ and $V$ which satisfy the standard commutation relation for the irrational rotation algebra 
$$UV\;=\;e^{2\pi i\delta}VU.$$  Rieffel's work showed that there was a unitary equivalence ${\mathcal W}:L^2(\R)\;\rightarrow\; L^2(\R)$ such that 
$${\mathcal W}U^{\delta}_{(1,0,0)}{\mathcal W}^{\ast}\;=\;V$$
and 
$${\mathcal W}U^{\delta}_{(0,1,0)}{\mathcal W}^{\ast}\;=\;U.$$

It then became clear that in order that there exist $f\in L^2(\R)$ such that 
$\{U^{\delta}_{(p,q,0)}(f): p,q\in \Z\}$ forms a NTF for $L^2(\R),$ it was necessary that the von Neumann algebra $M(\delta\Z, \R/\Z)$ should contain a cyclic vector in $L^2(\R).$ 
Now $M(\delta\Z, \R/\Z)$ and $M(\Z,(\delta\Z)\backslash \R)$ are both finite von Neumann algebras in the sense that they have defined on them canonical center-valued traces, which shall be denoted here by $\tau$ and $\tau',$ respectively.  Indeed, when $\delta$ is irrational, these von Neumann algebras are in fact $II_1$ factors, so that we can think of the traces as being scalar multiples of the identity.  Thus in the case of irrational $\delta,$ we can identify the range of these traces with complex numbers,  and in the case where the traces are defined on projections in the respective von Neumann algebras, the traces will take on real values between $0$ and $1.$  Fixing some notation, suppose that $M$ and its commutant $M'$ are two finite von Neumann algebras acting on a Hilbert space ${\mathcal H}.$ Then for $\xi\in {\mathcal H}$ let $E_{\xi}$ be that projection in the center of $M$ corresponding to the smallest closed subspace of ${\mathcal H}$ containing $M'(\xi),$ and let $E'_{\xi}$ be that projection in the center of $M'$ corresponding to the smallest closed subspace of ${\mathcal H}$ containing $M(\xi).$ We remark that if $\xi$ is a cyclic vector for $M,$ then $E'_{\xi}\;=\;\text{Id}_{M'}.$
Recall from \cite{Ri} that the coupling function $\omega$ is the possibly unbounded operator defined from the center of $M$ to the center of $M'$ which satisfies the condition
$$\tau'(E'_{\xi})\;=\;\tau(\omega(E_{\xi})).$$
Now suppose that $M(\delta\Z, \R/\Z)=M$ and $M(\Z,(\delta\Z)\backslash \R)=M'$
In Theorem 3.2 and Example 4.1 of \cite{Ri}, Rieffel computed the coupling function $\omega$ corresponding to the pair $(M, M')$ in this case to be exactly multiplication by the positive number $\frac{1}{\delta}.$ It follows that if $f$ is any vector in $L^2(\R),$
$$\tau'(E'_{f})\;=\;\frac{1}{\delta}\cdot\tau(E_f),$$
which is strictly less than one if $\delta>1$ (since $0\;\leq\;\tau(E_f)\;\leq 1$).
As $\tau'(E'_{f})<1,\;E'_f$ will not equal to $\text{Id}_{M'},$ which by our above remarks shows that $f$ cannot be cyclic for $M.$

In addition to the work of I. Daubechies, M. Rieffel and L. Baggett on the Gabor frame problem (sometimes called ``Weyl-Heisenberg frames" in the literature), Daubechies, H. Landau, and Z. Landau \cite{D3} also pursued different methods in dimension greater than $1,$ as did 
J. Ramanathan and T. Steger \cite{RS}.  Both of the papers cited above involved the use of von Neumann algebras in their proofs. A. Ron and Z. Shen \cite{RoS} also worked on this problem using a different approach.  
Very recently, D. Han and Y. Wang \cite{DL} have proved the following generalization of results on Gabor frames, which applies to lattices in $\R^n,$ for $n>1.$ 
\begin{theorem}
\label{thm hanwang}
\cite{HW}
Let $n\in \N,$ and let ${\mathcal L}$ and ${\mathcal K}$ be two full-rank lattices in $\R^n,$ called the modulation and translation lattices, respectively.
For $\lambda\in {\mathcal L},\kappa\in {\mathcal K},$ and $g\in\;L^2(\R^n),$ let 
$$g_{\lambda,\kappa}(x)\;=\;M_{\lambda}T_{\kappa}g(x),$$
where $M_{\lambda}$ and $T_{\kappa}$ are the standard modulation and translation operators.
Then there exists $g\in\;L^2(\R^n)$ such that 
$\{g_{\lambda,\kappa}:\lambda\in {\mathcal L}, \kappa\in {\mathcal K}\}$ form a frame for $L^2(\R^n)$ if
and only if $ v({\mathcal L})·v({\mathcal K})\leq 1.$
\end{theorem}
Here $v({\mathcal L})$ and $v({\mathcal K})$ denote  $|\det L|$ and $|\det K|$ respectively, where $L$ and $K$ are $n\times n$ matrices with ${\mathcal L}\;=\;L(\Z^n)$ and ${\mathcal K}\;=\;K(\Z^n).$ We remark that the necessity of the volume condition had been proven in \cite{RS} and \cite{RoS} and it is the sufficiency that was established in \cite{HW}. The proof of Han and Wang uses geometric lattice tiling and packing methods.

It would be interesting to approach the higher dimensional case by decomposing a certain representation of a discrete subgroup of the $2n+1$-dimensional Heisenberg Lie group as a direct integral.

\subsection{The study of wavelets using the direct integral of the wavelet representation}

Recall $G_A$ denotes the discrete semi-direct product group ${\Q}_A\rtimes_{\theta} \Z,$ where $A$ is an $n\times n$ integral dilation matrix and $\Q_A$ denotes the discretized $A$-adic subgroup of $\Q^n$ obtained from $A,$  viewed as a discrete group.  For example, if $n=1$ and $A$ is the constant dilation by $2,\;\Q_2$ represents the group of dyadic rational numbers.   

Recall from Definition \ref{def waverep} that the wavelet representation $W$ of $G_A$ is given by 
the formula 
$$W_{(\beta,k)}\;=\;T_{\beta}D_A^k,\;(\beta,k)\;\in\;G_A.$$
We discuss how to decompose this representation as a direct integral of representations.
The notion of wavelet set due to X. Dai and D. Larson in the classical $A=2$ case \cite{DL}, and for more  general dilation matrices to X. Dai, D. Larson and D. Speegle \cite{DLS} is needed.

We consider dilation and translation in the frequency domain, using the notation of Section 3.1.
Denote
$${\widehat{D_A}}={\mathcal F}D_A{\mathcal F}^{-1},\;\widehat{T_v}={\mathcal F}T_v{\mathcal F}^{-1},\;v\in\Z^n.$$
\begin{definition} 
\label{def waveset}
(\cite{DL, DLS})
Let $A$ be an $n\times n$ integer dilation matrix on $\R^n,$ with $A^T$ denoted by $B.$
We say that a set $E\;\subseteq\R^n$ of finite positive measure is a wavelet set for dilation by $A$ 
if $\frac{1}{\sqrt{\nu(E)}}\chi_E$ is a wavelet for dilation by $A$ in the frequency domain, that is, if 
$$\{{\widehat{D_A}}^j\widehat{T_v}(\frac{1}{\sqrt{\nu(E)}}\chi_E)\;:\;j\in\;\Z,\;v\in\;\Z^\}$$ is an orthonormal basis for $L^2(\R^n).$
\end{definition}
This definition would imply by use of the inverse Fourier transform that 
\newline ${\mathcal F}^{-1}(\frac{1}{\sqrt{\nu(E)}}\chi_E)$ is a wavelet for dilation by $A.$

Dai and Larson have shown that $E\;\subseteq\;\R^n$ is a wavelet set for dilation by $A$ if and only if  
the following conditions are satisfied:
\renewcommand{\labelenumi}{(\roman{enumi})}
\begin{enumerate}
\item  $E$ is a measurable subset of $\R^n$, with $B^j(E) \cap
B^k(E)=\emptyset$, $j \neq k \in \Z$,\item  $\nu(\R^n \backslash \cup_{j\in \Z} B^j(E))=0$;
\item  $E$ is translation congruent to the set $[-\frac{1}{2},\frac{1}{2})^n$ modulo the
lattice $\Z^n.$
\end{enumerate}
Conditions $(i)$ and $(ii)$ are equivalent to $E$ tiling  $R^n$ under dilation by $B,$ and condition $(iii)$ is equivalent to $E$ tiling the plane under translation by $\Z^n.$ 

The following theorem of X. Dai, D. Larson, and D. Speegle show that wavelet sets always exist.
\begin{theorem}
\label{thm wavesetexist}
\cite{DLS}
Let $A$ be an $n\times n$ integral dilation matrix.  Then there exists a wavelet set in $\R^n$ for dilation by $A.$
\end{theorem}
\begin{example}
\cite{DL}
Let $n=1,\;A=2$ and set $E=[-1,-1/2)\cup [1/2,1).$ It is clear that $E$ tiles the plane by dilation by $2$ and by translation by $\Z.$ Hence $E$ is a wavelet set and $\Psi(x)\;=\;\chi_E(x)$ 
 is a wavelet in the frequency domain for dilation by $2.$  The corresponding wavelet in the time domain is the Shannon wavelet related to the sinc function.
\end{example}

We now use wavelet sets and direct integral theory the wavelet representation.  We first consider a representation equivalent to the wavelet representation.
\begin{definition}
The representation ${\widehat{W}}: :\;{\Q}_A\rtimes \Z\rightarrow\;{\mathcal U}(L^2(\R^n))$  defined by 
$${\widehat{W}}(\beta,k)\;=\;{\mathcal F}W(\beta,k){\mathcal F}^{-1}$$ 
is called the {\bf wavelet representation of the wavelet group $G_A$ in the frequency domain}.
One calculates that 
$${\widehat{W}}(\beta,k)(f)(x)\;=\;|\text{det}(A)|^{-\frac{k}{2}}e^{2\pi i\beta x}f(B^k x),\;(\beta,k)\;\in\;G_A,\;f\in L^2(\R^n).$$
\end{definition} 

Now let $E$ be a wavelet set for translation by $A;$ by Theorem \ref{thm wavesetexist} such a set will exist.  Note that since $E$ tiles $\R^n$ by dilation by powers of $B,$ it is possible to identify $\R^n$ with the Cartesian product space $E\times \Z,$ and hence we can set up a unitary operator carrying $L^2(\R^n)$ to $L^2(E\times \Z)\cong L^2(E,l^2(\Z)),$ where $E$ carries ordinary Haar measure.  Upon doing this, and studying the wavelet representation of $G_A$ in the frequency domain, we obtain the following result:
\begin{lemma} (\ref{def waverep}, \cite{LPT})
\label{lem lptwavedecomp}
Let $A$ be a $n\times n$ integral dilation matrix and let $E$ be a wavelet set for dilation by $A.$ The wavelet representation $W$ of $G_A$ is equivalent to the representation ${\widetilde{W}}:\;G_A\;\rightarrow\;{\mathcal U}(L^2(E\times\Z))$
defined by 
$$\widetilde{W}_{(\beta,k)}(f)(x,m) = e^{-2\pi i x\cdot A^m(\beta)}f(x,m-k),\;(\beta,m) \in \Q_A,\;f\in L^2(E\times \Z).$$
\end{lemma}
The lemma above allows us to decompose the wavelet representation of $G_A$ as a direct integral of monomial representations, just as Baggett did earlier in the Gabor case.  The following result is joint with L.-H. Lim and K. Taylor \cite{LPT}:
\begin{theorem}
\label{thm waverepdecomp}
\cite{LPT}
Let $A$ be a $n\times n$ integral dilation matrix, and let $E$ be a wavelet set for dilation by $A$ in the sense of \ref{def waveset}.  Then:
\renewcommand{\labelenumi}{(\roman{enumi})}
\begin{enumerate}
\item The wavelet representation  $W$ of $G_A$ is equivalent to a direct integral of representations $\int^{\oplus}_E{\widetilde{W_x}} d\nu$ of $G_A$ on the direct integral Hilbert space $\int^{\oplus}_E [l^2(\Z)] d\nu.$ 
\item Each ${\widetilde{W}}$ is a representation of $G_A$ on $l^2(\Z)$ which is induced from a character $\chi_x$ on the subgroup $\Q_A$ of $G.$ If $x\not= 0,\;{\widetilde{W_x}}$ is irreducible.
\end{enumerate}
\end{theorem}
We briefly discuss the proof of the above theorem.  Firstly, the discussion on direct integrals in Section $2$ shows that the Hilbert space $L^2(E\times\Z)$ can be identified with $L^2(E,l^2(\Z))$ which in turn can be identified with the direct integral Hilbert space $\int_E^{\oplus}[l^2(\Z)]d\nu.$  Upon making this identification, it becomes clear that the 
 representation ${\widetilde{W}}$ can be decomposed as a direct integral in the desired fashion.
The fact that each ${\widetilde{W_x}}$ is induced from a character on $\Q_A$ follows from calculations done in \cite{LPT}. The irreducibility of the representation ${\widetilde{W_x}}$ for $x\not= 0$ follows from the observation that the orbit of $\chi_x$ under the induced action of $\Z$ in the compact abelian dual group $\widehat{\Q_A}$ is dense whenever $x\not= 0.$

Before making some more observations, we recall the definition of the (left) regular representation of a group $G.$
\begin{definition}
Let $G$ be a second countable locally compact group.  The (left) regular representation of $G$ is the representation $\lambda$ of $G$ on $L^2(G)$ given by 
$$\lambda_g(f)(h)\;=\;f(g^{-1}h),\;g,h\;\in G,\;f\;\in\;L^2(G).$$
(Here G is given its left-invariant Haar measure).
\end{definition}

We remark that in their original study of wavelet sets \cite{DL}, Dai and Larson actually computed the commutant of the wavelet representation of $\Q_2\rtimes\Z$ and showed that it was a commutative von Neumann algebra.  They also noted that the von Neumann algebra generated by the regular representation of $\Q_2\rtimes\Z$ was a $II_1$-factor, as was its commutant.  Hence the wavelet representation and the regular representation of $\Q_2\rtimes\Z$ cannot be unitarily equivalent. The same argument shows that similar statements are true for the regular representation and the wavelet representation of $\Q_A\rtimes_{\theta}\Z.$ By our discussion in Section $2,$ this implies that there is up to unitary equivalence, only one direct integral decomposition of the wavelet representation $W$ of $G_A,$ whereas using the method outlined in \cite{Mac}, one can decompose the regular representation of $G_A$ into two distinct direct integral decompositions which are not equivalent. 

One can use Lemma \ref{thm waverepdecomp} to prove the following Corollary, which in the classical case is due to F. Martin and A. Vallette \cite{MV}. We first recall that two unitary representations of a locally compact group $G$ are said to be {\bf weakly equivalent} if the kernels of the associated representations of the group $C^{\ast}$-algebra $C^{\ast}(G)$ are equal. 
\begin{corollary} 
\cite{LPT}
The wavelet representation of $G_A$ is weakly equivalent to the regular representation of $G_A.$
\end{corollary}

We end this section by discussing an example which shows how to construct tight frames in the classical case, when $n=1$ and $A$ is dilation by $2,$ by using the direct integral approach.  This example is analogous to the approach of Baggett in the Gabor frame case.

\begin{example}
\label{ex waveframe} 
Fix $a,b>0$ and let $E= [-2a,-a)\cup[b,2b)\subseteq \R$. Let $A=(2);$ 
note that dilation of $E$ by $2$ tiles the real line
$\R$. We claim that $\psi = \mathcal{F}^{-1}(1/\sqrt{\mu(E)}\chi_E)$ is a
tight frame wavelet for translation by the integers and dilation by 2 if and only
if $a+b\leq 1$.  If $a+b< 1$, the vectors $\{ D_2^m T_k\psi \mid
m,k \in \Z)\}$ form a tight frame for $L^2(\R)$. If $a+b=1\;\psi$ is an orthonormal wavelet for dilation by $2.$

We do not do the calculations in detail, but note the main idea for this calculation is very similar to that given in Example \ref{ex Gabor} for the Gabor frame case.
We begin with the case $a+b=1$; then setting $b=1 - a$, it is 
already known (c.f. \cite{DL}, Example 4.5) that $[-2a,-a)\cup [1 -a,2 -2a)$ is
a wavelet set for dilation by $2,$ and the statement follows.   If $a+b < 1$,
$E$ still tiles $\R$ under dilation by $2$. Using the methods of Lemma \ref{lem lptwavedecomp}, it is possible to show that the wavelet representation is unitarily equivalent to a representation over a direct integral Hilbert space 
$\int^{\oplus}_{E}[l^2(\Z)]d\nu.$  Using this representation, the functions   
 $\{D_{2}^m T_k\psi:\;m,k \in \Z \}$ correspond to the functions $\{\frac{1}{\sqrt{\mu(E)}}e^{-2\pi ikx}\otimes\delta_m:\; 
k\in \Z \}\;\subseteq\;L^2(E)\otimes l^2(\Z)\cong\;\int^{\oplus}_{E}[l^2(\Z)]d\nu.$  Since $\{\delta_m :\; m\in \Z\}$ is an orthonormal basis for $l^2(\Z),$ it follows that $\{D_{2}^m T_k\psi;\; m,k \in \Z \}$ will form a frame for
$L^2(\R)$ if and only if the functions $\{\frac{1}{\sqrt{\nu(E)}}e^{-2\pi ikx};\;
k\in \Z \}$ form a frame for $L^2(E)$. If $a+b < 1$,
one computes that $\{\frac{1}{\sqrt{\nu(E)}}e^{-2\pi ikx};\;
k\in \Z \}$ form a tight frame for $L^2(E)$ with frame constant $\nu(E)=a+b$ just as in Example 
\ref{ex Gabor}. 

Finally, if $a+b =\nu(E)>1 $, the closed span of the set $\{1/\sqrt{\mu(E)}e^{-2\pi ikx}
\mid k\in \Z \}$ in $L^2(E)$ will consist of those $L^2$ functions on $E$
which are $\Z$-periodic. Hence this closed span cannot be equal to
$L^2(E)$, so that the functions $\{1/\sqrt{\mu(E)}e^{-2\pi ikx} \mid k\in \Z
\}$ do not form a frame for $L^2(E).$ Just as in Example \ref{ex Gabor}, this implies that the set $\{D_{(2)}^m T_k\psi \mid m,k\in \Z \}$ cannot form a frame for $L^2(\R)$ in this case.
\end{example}
    

\section{Multiresolution Analysis, Generalized Multiresolution Analysis, and the SNAG Theorem}

We now discuss multiresolution analyses and generalized multiresolution analyses. The concept of multiresolution analysis (commonly abbreviated by the acronym MRA) was originally developed by S. Mallat and Y. Meyer in the mid- 1980's (\cite{Ma,Me} as a new way to construct interesting frames and wavelets.   {\it A priori}, it would appear that their construction is not connected to the work of Stone and/or von Neumann.  If one looks at the basic algorithm from a slightly different approach, especially in the generalized multiresolution analysis case, it becomes clear that there is a connection.  We first review the Stone-Naimark-Ambrose-Godement generalization of the Spectral Theorem for bounded self-adjoint operators due to D. Hilbert, E. Hellinger, and H. Hahn. 


\subsection{Unitary representations of locally compact abelian groups: the SNAG theorem}

In 1930, building on earlier observations of H. Weyl, M. Stone noted the connection between unitary representations of $\R$ on a Hilbert space ${\mathcal H},$ unbounded self-adjoint operators on ${\mathcal H},$ and projection valued measures (resolutions of the identity) on ${\mathcal H}.$ To be more specific, he proved that associated to continuous unitary representation $t\mapsto U_t$ of $\R$ on ${\mathcal H}$ there is a unique (possibly unbounded) self-adjoint operator $A$ defined on ${\mathcal H}$ such that for every $v,w\in{\mathcal H},\forall t\in\R,$
$$<U_t(v),W>\;=\;e^{2\pi i t<Av,w>}.$$ We express this in operator terms by writing 
$$U_t\;=\;e^{2\pi i tA},\;t\in A.$$
By the ordinary spectral Theorem for unbounded self-adjoint operators, we can write 
$$A\;=\;\int_{\R}xdP_{x},$$ where $P$ is a projection valued measure on $(\R,{\mathcal B}),$ where ${\mathcal B}$ denotes the Borel subsets of $\R.$
Hence we can write 
$$U_t\;=\;\int_{\R}e^{2\pi itx}dP_{x},\;\forall t\in\R.$$
By results relating projection valued measures to direct integral theory (see Theorem III.13 of Fabec's book \cite{Fa}) this implies that there exists 
a Borel measure $\nu$ on $(\R,{\mathcal B})$ such that 
${\mathcal H}$ can be decomposed as a direct integral 
$${\mathcal H}\;=\;\int^{\oplus}_{\R}[{\mathcal H}_x] d\nu.$$ 
The class of the measure $\nu$ is unique.  With respect to this direct integral decomposition of ${\mathcal H},$ we can write
\begin{equation}
\label{eq Stonethm}
U_t\;=\;\int_{\R}^{\oplus}[e^{2\pi itx}\;\text{Id}]d\nu(x),\;t\in\R,
\end{equation}
where $\text{Id}$ denotes the identity operator on ${\mathcal H}_x.$
Now let $\hat{G}$ denote the Pontryagin dual group of characters of a locally compact abelian group $G.$ Note that because $\R\;\cong\;\widehat{\R}$ via the correspondence $x\mapsto\;e^{2\pi ix\cdot}\;=\chi_x,$ Equation \ref{eq Stonethm} can be rewritten more generally as 
$$U_t\;=\;\int_{\widehat{\R}}^{\oplus}[\chi(t)\;\text{Id}]d\nu(\chi),\;\chi\in\widehat{\R},\;t\in\R,$$ for a suitable Borel measure on $\widehat{\R}.$

This result was generalized to an arbitrary locally compact abelian groups by W. Ambrose \cite{A}, M. Naimark \cite{N} and R. Godement \cite{G} in the early 1940's. Indeed G. Mackey calls the general result the ``SNAG" theorem, for obvious acronymic reasons.

\begin{theorem}
\label{thm SNAG}
\cite{S,N,A,G}
Let $G$ be a locally compact abelian group, and let $U$ be a continuous unitary representation of $G$ on the Hilbert space ${\mathcal H}.$  Let $\hat{G}$ denote the Pontryagin dual group of $G.$  Then there is a Borel measure $\nu$ on $\hat{G}$ such that ${\mathcal H}$ can be decomposed as a direct integral 
$${\mathcal H}\;=\;\int^{\oplus}_{\hat{G}}[{\mathcal H}_{\chi}] d\nu(\chi),$$ 
and with respect to this decomposition we can write 
$$U_g\;=\;\int_{\hat{G}}^{\oplus}[\chi(g)\text{Id}]d\nu(\chi),\;\chi\in\hat{G},\;g\in G.$$
The measure class of $\nu$ is unique.
\end{theorem}
This result specialized to the case $G\;=\; \R$ gives Stone's Theorem \cite{S}.

In wavelet theory, this result will most often be applied to the case where $G\;=\;\Z^n,$ and ${\mathcal H}$ is a closed subspace of $L^2(\R^n),$ usually denoted by $V_0,$ which is invariant under the translation operators $\{T_v:\;v\in \Z^n\}.$


\subsection{Multiresolution Analysis, and the construction of wavelet and frame bases}

We give now the definition of Multiresolution Analysis due to S. Mallat and Y. Meyer. 
\begin{definition}
\label{def multires}  
 Let $\{V_i\}_{i\in\Z}$ be a bisequence of closed subspace of $L^2(\R^n),$ and let $A$ be an $n\times n$ dilation matrix, that is, a matrix all of whose eigenvalues have modulus greater than $1.$ We say that 
$\{V_i\}_{i\in\Z}$  is {\bf multiresolution analysis (MRA)} for dilation by $A$ if:
\renewcommand{\labelenumi}{(\roman{enumi})}
\begin{enumerate}
\item $\;\;\cdots V_{-1}\subseteq V_0\subseteq V_1\cdots$ (the $V_i$ form a nested sequence of closed subspaces of $L^2(\R^n));$
\item $D_A^i(V_0)\;=\;V_i,\;i\in\Z;$
\item $\overline{\cup_{i\in\Z}V_i}\;=\;L^2(\R^n),\;\cap_{i\in\Z}V_i=\{0\};$
\item There exists $\phi\in L^2(\R),$ called a scaling function for dilation by $A,$ such that
\nn $\{T_v(\phi):\;v\in\Z^n\}$ is an orthonormal basis for $V_0.$
\end{enumerate}
\end{definition}
The last condition can be interpreted as meaning that $V_0$ can be identified with $l^2(\mathbb Z^n),$ and hence by Pontryagin duality, with 
$L^2(\mathbb T^n).$  Using this identification, we see that the representation of $\mathbb Z^n$ on $\widehat{V_0}$ given by the operators $\{\widehat{T_v}\}$ is equivalent to the regular representation of $\Z^n$ on $l^2(\Z^n).$
From the direct integral point of view, using the SNAG theorem just discussed, this is equivalent to saying that  the scaling function gives us a spectral decomposition of the translation operators $T_v,\;v\in\Z^n,$ acting on  of the space $V_0.$  Let us set $n=1$ to illuminate this idea more clearly.  The measure on $\hat{\Z}\cong\T$ is guaranteed by the SNAG theorem turns out to be the Haar measure on $\T.$ Define a map 
$J:\;V_0\;\rightarrow\; L^2(\mathbb T)$ by 
$$J(\sum_{k\in\mathbb Z}c_k T^k(\phi))(x)\;=\;\sum_{k\in\mathbb Z}c_ke^{2\pi ikx}.$$
Then using elementary facts about orthonormal bases in Hilbert spaces, it is easy  to see $J$ is a unitary isomorphism, and that for every $f\in V_0,$
$$J(T(f))(x)\;=\;e^{2\pi ix}J(f)(x).$$
We will return to this construction when discussing the generalized multiresolution analyses of L. Baggett, H. Medina, and K. Merrill \cite{BMM}. 

We first view MRA's in the frequency domain by means of the Fourier transform.
\begin{definition}
\label{def MRAfreq}
Let $\{\widehat{V_i}\}_{i\in\Z}$ be a bisequence of closed subspace of $L^2(\R^n).$ 
We say that $\{\widehat{V_i}\}_{i\in\Z}$  is a {\bf multiresolution analysis in the frequency domain} for dilation by $A$ if:
\renewcommand{\labelenumi}{(\roman{enumi})}
\begin{enumerate}
\item $\;\;\cdots \widehat{V_{-1}}\subseteq\widehat{V_0}\subseteq\widehat{V_1}\cdots$ (the $\widehat{V_i}$ form a nested sequence of closed subspaces of $L^2(\R^n));$
\item $(\widehat{D_A})^i(\widehat{V_0})\;=\;\widehat{V_i},\;i\in\Z;$
\item $\overline{\cup_{i\in\Z}\widehat{V_i}}\;=\;L^2(\R^n),\;\cap_{i\in\Z}\widehat{V_i}=\{0\};$
\item There exists $\Phi\in L^2(\R^n),$ a scaling function in frequency domain, such that
\nn $\{\widehat{T_v}(\Phi):\;v\in\Z^n\},$ is an orthonormal basis for $\widehat{V_0}.$
\end{enumerate}
\end{definition}
Of course one moves between MRA's  and MRA's in the frequency domain by setting 
$$\widehat{V_i}\;=\;{\mathcal F}(V_i),\; V_i\;=\;{\mathcal F}^{-1}(\widehat{V_i})$$
\begin{example} 
Let $n=1$ and $A=(2).$ Then
$$\Phi(x)\;=\;\left\{\begin{array}{ll}
\frac{e^{2\pi ix}-1}{2\pi ix},&\mbox{if}\ \;x\not= 0,\\
1&\;\mbox{if}\ \;x\;=\;0
\end{array}\right.$$
is a scaling function in the frequency domain for dilation by $2.$
\end{example}
Let $A$ be a $n\times n$ integral dilation matrix. Given a MRA for dilation by $A,$ there is an algorithm using the scaling function $\Phi$ for constructing a wavelet family
in the frequency domain for dilation by $N.$
In this algorithm, the scaling function $\Phi$ is extremely important, since from $\Phi$ one constructs $\widehat{V_0}$ using the translation operators $\widehat{T_v},$ and then the other subspaces $\widehat{V_i}$ using the dilation operator $\widehat{D},$ and finally the wavelet family from looking at $\widehat{W_0}=\widehat{V_1}\ominus\widehat{V_0}.$  

In order to construct $\Phi,$ it helps to have a {\bf low-pass filter} $m_0$ for dilation by $A.$ For simplicity of exposition we let $n=1$ and let the dilation matrix be $N$ where $N$ is a positive integer greater than $1.$

\begin{definition}
\label{def lo-pass}
 Let $N$ be a positive integer greater than $1.$  A {\bf low-pass filter} $m_0$ for dilation by $N$ is a $\Z$-periodic function $m_0:\;\R^n\;\rightarrow\;\C$ which satisfies the following conditions:
\renewcommand{\labelenumi}{(\roman{enumi})}
\begin{enumerate}
\item $m_0(0)\;=\;\sqrt{N};$
\item $\sum_{l=0}^{N-1}|m_0(x+l/N)|^2=N$ a.e.
\item $m_0$ is Holder continuous at $0$ and is non-zero in a sufficiently large neighborhood of $0$ (``Cohen's condition", c.f. \cite{Co,D2})
\end{enumerate}
\end{definition}
Given a low-pass filter, one can use it to construct the scaling function in the frequency domain as follows:
\begin{theorem} 
\cite{Ma,Me,D2}
 Given a low-pass filter $m_0$ for dilation by $N$ which satisfies Cohen's condition, then
$$\Phi(x)\;=\Pi_{i=1}^{\infty}[\frac{m_0(N^{-i}(x))}{\sqrt{N}}]$$
converges a.e. and in $L^2(\R)$ to a scaling function in the frequency domain for dilation by $N,$ which can be used to construct an orthonormal wavelet family for dilation by $N.$
\end{theorem}

We briefly discuss how one uses the scaling function, the low-pass filter, and other filters called high-pass filters, to construct the wavelet family. 
\begin{definition}  
\label{def hipass}
Let $n=1$ and let the dilation matrix $A$ now be the positive integer $N,$ where $N>1,$ and let  $m_0$  be a  low-pass filter for dilation by $N.$ A set of essentially bounded  measurable $\Z$ periodic functions $m_1,m_2,\cdots,m_{N-1}$ are called {\bf high-pass filters} associated to $m_0,$ if the 
$N\times N$ matrix
$$(\frac{m_i(x+j/N)}{\sqrt{N}})_{0\leq i,j \leq N-1}\;\in\;{\mathcal U}(N, L^{\infty}(\R/\Z)).$$
\end{definition}

The following result is due to Y. Meyer and S. Mallat, and in the normalized tight frame case to W. Lawton.
\begin{theorem}
\cite{Ma,Me,L}  Let $N$ be a positive integer greater than $1,$ let $m_0$ be a  low-pass filter  for dilation by $N,$ and let 
$\Phi$ be the scaling function constructed from $m_0$ above. If  $m_1,m_2,\cdots,m_{N-1}$ are high-pass filters associated to $m_0,$ then 
$$\{ \Psi_1=\widehat{D}(m_1\Phi),\;\Psi_2=\widehat{D}(m_2\Phi),\;\cdots,\;\Psi_{N-1}=\widehat{D}(m_{N-1}\Phi)\}$$ is a orthonormal wavelet family in the frequency domain for dilation by $N.$ If Cohen's condition is not satisfied, the $\{\Psi_k\}$'s still form a NTF wavelet family. 
\end{theorem}

The above discussion indicates how to obtain MRA's from filter functions.
It may not be immediately apparent how filter functions are related to the work of von Neumann, but a very nice paper by L. Baggett, A. Carey, W. Moran, and P. Orring \cite{BCMO} gives an exposition of wavelet theory from the point of view of von Neumann algebras, and shows how one can obtain filter functions from MRA's by using a result on cancellation for finite von Neumann algebras. The result of Baggett, Carey, Moran and Orring employed von Neumann algebras to study certain lattices $\Gamma$ of locally compact abelian groups $G$ having a faithful unitary representation ${\mathcal U}$ on a Hilbert space ${\mathcal H}$ via translation. They also assumed the existence of a generalized ``dilation operator" $\delta$ acting on ${\mathcal H}$ with the property that 
$\{\delta^{-1}{\mathcal U}_{\gamma}\delta:\; \gamma\;\in\;\Gamma\}$ was a subgroup of finite index in $\{{\mathcal U}_{\gamma}:\; \gamma\;\in\;\Gamma\}.$ Here for simplicity $\Gamma=\Z,\;G=\R,\;{\mathcal H}\;=\;L^2(\R),$ and $\Z$ acts on $L^2(\R)$ by translation operators. We also assume that the dilation operator $\delta$ is the standard dilation $D_N$ corresponding to the positive integer $N>1,$ which we denote here by $D.$ We refer the reader to \cite{BCMO} for the general set-up; we only note that the main theme of the argument carries over to the general case.  

Let $\{V_j\}_{j\in\Z}\;\subseteq L^2(\R)$ be an ordinary multiresolution analysis for dilation by $N.$  Then $V_0$ is invariant under translations $\{T_v: v\in\Z\},$ and by definition, the corresponding representation of $\Z$ on $V_0$ is equivalent to the regular representation of $\Z$ on $l^2(\Z).$  Let $W^{\ast}(\Z)$ denote the commutative von Neumann subalgebra of ${\mathcal B}(L^2(\R))$ generated by $\{T_v : v\in\Z\}.$ Then $V_1\;=\;D(V_0)$ is invariant under the action of $W^{\ast}(\Z)$ with invariant subspace $V_0,$ and the multiresolution theory tells us that 
$$V_1\;=\;V_0\oplus\;W_0,$$ where $W_0$ is the wavelet space, which is an invariant subspace for the von Neumann algebra $W^{\ast}(\Z).$  Using cancellation properties for finite von Neumann algebras (c.f. \cite{Dx}, Section 3.2.3 Proposition 6), they are able to show the following result:  
\begin{theorem}
\label{thm bcmo}
\cite{BCMO}
The wavelet space $W_0$ can be decomposed into $N-1$ invariant orthogonal subspaces $\{W_0^i\}_{i=1}^{N-1},$ each invariant under the action of $W^{\ast}(\Z).$  For each $i,$ the action of $\Z$ on the representation of $\Z$ on $W_0^i$ corresponding to the action of $W^{\ast}(\Z)$ is equivalent to the regular representation of $\Z$ on $l^2(\Z).$
\end{theorem} 

By definition $V_0$ contains an element $\phi$ such that $\{T_v(\phi)\}_{v\in\Z}$ forms an orthonormal basis for $V_0.$
Viewing this in the Fourier domain,
$$\{\widehat{T_v}(\widehat{\phi})\}_{v\in\Z}\;=\;\{\widehat{T_v}(\Phi)\}_{v\in\Z}$$ forms an orthonormal basis for $\widehat{V_0}={\mathcal F}(V_0),$ where $\Phi\;=\;\widehat{\phi}.$
Since $\Phi\;\in\;\widehat{V_1}\;=\;\widehat{D}(\widehat{V_0}),$ we can write 
$$\Phi(x)\;=\;\;\widehat{D}[{\mathcal F}(\sum_{v\in\Z}a_{0,v}T_v(\phi))](x)\;=\;\widehat{D}[(\sum_{v\in\Z}a_{0,v}\;e^{2\pi ivx})\Phi](x),$$
where the infinite sums converge in the relevant Hilbert space norm.
Writing $$m_0(x)\;=\;\sum_{v\in\Z}a_{0,v}\;e^{2\pi ivx},$$
we get $$\Phi(x)\;=\;\widehat{D}(m_0\Phi)(x)\;=\;\frac{1}{\sqrt{N}}m_0(\frac{x}{N})\Phi(\frac{x}{N}),$$
which is exactly the refinement equation in the frequency domain.
The function $m_0$ turns out to be the low-pass filter in this situation.
 
Similarly, each subspace $W_0^i,\;1\leq i\leq N-1,$ contains an element $\psi_i$ such that 
$\{T_v(\psi_i)\}_{v\in\Z}$ forms an orthonormal basis for $W_0^i.$
Again we have  $$\{\widehat{T_v}(\widehat{\psi_i})\}_{v\in\Z}\;=\;\{\widehat{T_v}(\Psi_i)\}_{v\in\Z}$$ forming  an orthonormal basis for $\widehat{W_0^i}={\mathcal F}(W_0^i),$ where $\Psi_i\;=\;\widehat{\psi_i}.$
Since $\Psi_i\;\in\;\widehat{V_1}\;=\;\widehat{D}(\widehat{V_0}),\;i\in\{1,\cdots N-1\},$ 
we see that $\xi_i\in\widehat{V_0}$ with $\psi_i\;=\;\widehat{D}(\xi_i).$
Then as before,  
$$\xi_i(x)\;=\;{\mathcal F}(\sum_{v\in\Z}a_{i,v}T_v(\phi))\;=\;[\sum_{v\in\Z}a_{i,v}\;e^{2\pi ivx}]\Phi(x).$$
Writing $$m_i(x)\;=\;\sum_{v\in\Z}a_{i,v}\;e^{2\pi ivx},\;1\;\leq\;N-1,$$
we obtain $$\Psi_i(x)\;=\;\widehat{D}(m_i\Phi)(x),$$ and the 
family $\{m_1,\;m_2,\;\cdots, m_{N-1}\}$ are 
the high-pass filter functions associated to $m_0.$
Standard arguments from wavelet theory (c.f. \cite{Str}) show that the filter functions $\{m_0,\;m_1,\;m_2,\;\cdots, m_{N-1}\}$ constructed in this fashion satisfy the standard high-pass orthogonality conditions as condensed in Definitions \ref{def lo-pass} and \ref{def hipass}.

An approach to finding continuous high-pass filter functions corresponding to a given low-pass filter function by using a $C^{\ast}$-module point of view, which has features in common with the approach in \cite{BCMO}, can be found in a recent article by the author and M. Rieffel \cite{pr}.

\subsection{Generalized multiresolution analyses and relationships to the SNAG Theorem, and current research}

In the frequency domain version of MRA for $n=1,$ Condition (iv) implies that the representation of $\Z$ on 
the initial space $\widehat{V_0}$ is equivalent to the regular representation of $\Z$ on $l^2(\Z).$
This means that when we use the Spectral Theorem to rewrite this representation as a direct integral on $\widehat{V_0},\;$ there is a measure $d\nu$ on $\hat{\Z}\;\cong\;\T\;\cong\;\R/\Z,$ which we parametrize by $x\in\;[0,1),$ such that 
$$\widehat{V_0}\;\equiv\;\int^{\oplus}_{\T}[\C]\;d\nu(x)$$ and 
$\widehat{T_v}\;=\;\int^{\oplus}_{\T}M_{e^{2\pi ivx}}\;d\nu(x).$ It is clear from our Fourier transform of $T$ that $d\nu(x)$ is just Haar measure $dx$ on $\T.$
The use of the SNAG Theorem comes when the wavelet theory is extended to translations $T_v,\;v\;\in\;\Z^n,$ acting on $L^2(\R^n),$ and more general dilations $D_A$ corresponding to $n\times n$ integer dilation matrices $A.$ In the ordinary MRA case, we get a somewhat analogous situation to the case when $n=1.$ Things become more complicated when condition (iv) of the definition of an MRA is not satisfied.

We give an example of a wavelet which does not come from an MRA.
\begin{example} 
\label{ex Journe}
Consider the wavelet in the frequency domain coming from the wavelet set 
$$E\;=\;[-\frac{16}{7},-2)\cup [-\frac{1}{2},-\frac{2}{7})\cup [\frac{2}{7},\frac{1}{2})\cup [2,\frac{16}{7}).$$
Then $\chi_E=\Psi_J$ is known as the {\bf Journ\'e wavelet}, and  is a wavelet for dilation by $2$ in the frequency domain that does not come from an MRA .  However, one can still form the nested sequences of Hilbert spaces 
$$V_i\;=\;{\overline{\text{span}}}\{D^jT_v({\mathcal F}^{-1}(\Psi_J))\;|\;v\in\Z,\;j\leq i\}$$
and it is possible to show:
\renewcommand{\labelenumi}{(\roman{enumi})}
\begin{enumerate}
\label{gmraconds}
\item $\;\;\cdots V_{-1}\subseteq V_0\subseteq V_1\cdots$ (the $V_i$ form a nested sequence of closed subspaces of $L^2(\R));$
\item $D^i(V_0)\;=\;V_i,\;i\in\Z;$
\item $\overline{\cup_{i\in\Z}V_i}\;=\;L^2(\R),\;\cap_{i\in\Z}V_i=\{0\};$
\item $V_0$ is invariant under all powers of $T.$
\end{enumerate}
\end{example}

The situation when one has a wavelet family and NTF wavelet frame which does not come from an MRA has been approached in various ways in independent works of  L. Baggett, H. Medina and K. Merrill \cite{BMM}, J. Benedetto and S. Li \cite{BL}, and D. Han, D Larson, M. Papadakis, and T. Stavropoulos \cite{HLPS}. We discuss here the approach due Baggett, Medina and Merrill. 
\begin{definition}
\label{def GMRA2}
Let $A$ be a $n\times n$ integral dilation matrix. A {\bf generalized multiresolution analysis} (GMRA) for dilation by  $A$ is a sequence of closed subspaces $\{V_i\}_{i\in\Z}$ of $L^2(\R^n)$ which satisfy the following conditions:
\renewcommand{\labelenumi}{(\roman{enumi})}
\begin{enumerate}
\item $V_i\subseteq V_{i+1}$ for all $i\in\Z;$
\item $D_A(V_i)\;=\;V_{i+1},\;\forall i\in\Z;$
\item $\overline{\cup_{i\in\Z}V_i}\;=\;L^2(\R^n),\;\cap_{i\in\Z}V_i=\{0\};$
\item $V_0$ is invariant under all of the operators $\{T_v:\;v\in\Z^n\}.$
\end{enumerate}
\end{definition} 
One can reformulate the above conditions in the frequency domain if desired, just as in Definition \ref{def MRAfreq}.

Baggett, Merrill and Medina showed given a wavelet family, one can correspond to it a GMRA, as in Example \ref{ex Journe}.  Given a GMRA not associated to a wavelet, they also 
were able to develop several characteristic invariants associated to a GMRA using spectral theory. They generalized the unitary operator $J$ defined from $V_0$ to $L^2(\mathbb T)$ in the classical MRA case to the GMRA case in the following theorem.
\begin{theorem} 
\cite{BMM} 
Given a GMRA in $L^2(\mathbb R^n)$ corresponding to dilation $D_A$ by an integral dilation matrix $A$ and the transform of integer translation $T_v,$ there is a unique sequence of Borel subsets $S_1\supseteq S_2\supseteq\cdots$ of $\mathbb T^n$ and a unitary operator $J:V_0\;\rightarrow\;\oplus_j\;L^2(S_j)$ such that 
$$[J(T_v(f)]_j(x)\;=\;e^{2\pi ix\cdot v}[J(f)]_j(x),\;f\in L^2(\T^n),\;v\;\in\;\Z^n.$$
\end{theorem}
We briefly discuss the construction of $J$ and the sets $\{S_j\}$ from the point of view of the SNAG Theorem. 
The representation $T$ of $\Z^n$ on $V_0$ can be decomposed as a direct integral.  With respect to this direct integral we can write $V_0\;\cong\;\int^{\oplus}_{\widehat{\Z^n}}[{\mathcal H}_x]d\nu(x),$
where $\nu$ is a Borel measure on $\widehat{\Z^n}=\T^n,$ and 
then by the SNAG theorem, with respect to this decomposition we can write 
$$T_v\;=\;\int^{\oplus}_{\mathbb T^n}e^{2\pi ix\cdot v}\text{Id}]d\nu(x),\;v\;\in\;\Z^n.$$
Baggett Merrill and Medina first proved that the class of the measure $\nu$ coming from the SNAG theorem is exactly the class of 
Haar measure on $\T^n,$ and then set 
$$S_j\;=\;\{x\in\;\T^n\;\text{dim}({\mathcal H}_x)\;\geq\;j\}.$$ 

Baggett, Merrill and Medina called the function $\mu(x)\;=\;\sum_j\chi_{S_j}(x)$ defined on $\mathbb T$ the {\bf multiplicity function} corresponding to the GMRA $\{V_i\}_{i\in\Z},$
because $\mu(x)$ indicates the multiplicity of the character 
$e^{2\pi ix\cdot}\in\widehat{\Z^n}=\T^n$ in the above decomposition of the representation $\T^n,$ i.e.,
$$\mu(x)\;=\;\text{dim}({\mathcal H}_x),\; x\in\T^n.$$

For simplicity, we specialize to the case where $n=1$ and $A$ is dilation by a positive integer $N>1.$  Baggett, Merrill and Medina showed that $\mu$ satisfies 
$$\mu(x)\;\leq\;\sum_{l=0}^{N-1}\mu(\frac{x+l}{N})\;\text{a.e.}.$$
 If $\mu$ is essentially bounded, with $c=\text{ess sup}\mu(x),$ it is possible to define the conjugate multiplicity function 
$$\tilde{\mu}(x)\;=\;\sum_{l=0}^{N-1}\mu(\frac{x+l}{N})-\mu(x)$$
By definition, $\mu$ and $\tilde{\mu}$ satisfy the so-called ``consistency equation":
\begin{equation}
\label{coneq}
\mu(x)\;+\;\tilde{\mu}(x)\;=\;\sum_{l=0}^{N-1}\mu(\frac{x+l}{N}).
\end{equation}

Baggett, J. Courter and Merrill in \cite{BCM} then generalized the Mallat and Meyer algorithm for constructing wavelets from filters to this GMRA setting.  They first generalized the concept of low and high pass filters.  Given an integer valued function $\mu$ (essentially bounded by c on $\mathbb T$ and satisfying technical conditions outlined in 
\cite{BM} they defined ``generalized conjugate mirror filters", which correspond to low-pass filters in the classical case, to be functions $\{h_{i,j}\}_{1\leq\;i,j\;\leq c},$ where each $h_{i,j}$ is supported on $S_j.$ 
Similarly, they defined ``complementary conjugate mirror filters," an analogue of classical high pass filters, to be functions
$\{ g_{k,j}\}_{1\leq\;k\;\leq d,\;1\leq\;j\;\leq c},$ where each $g_{k,j}$ is supported on $S_j.$ 

The functions $\{h_{i,j}\}$ and $\{ g_{k,j}\}$ satisfy orthogonality conditions which are modified versions of the orthogonality conditions satisfied by classical low and high-pass filters; see \cite{BCM} for precise details. They built examples of functions $g_{k,j}$ and $h_{i,j}$ satisfying these conditions by using an explicit algorithm, and then modified the results obtained  to get examples of generalized filters, mimicking the way examples of filters are obtained in the classical case \cite{BCM}.  

Under appropriate conditions on these generalized filter functions, Baggett, Courter and Merrill then used them to construct a finite tight frame wavelet family $\{\Psi_1,\cdots,\Psi_d\}\subseteq L^2(\mathbb R)$. 

They first built generalized scaling functions 
$\{\Phi_1,\Phi_2, \cdots, \Phi_c\}\;\subseteq\;\widehat{V_0}$ using an infinite product construction involving dilates of a matrix with periodizations of $\{h_{i_j}\}$ as entries (c.f. Theorem 3.4 \cite{BCM}). The $\{\Phi_i\}_{i=1}^c$ appear as the entries in the first column in the infinite product matrix. 

Given the above notation and construction, they then have:
\begin{theorem} \cite{BCM} Let  $\{(h_{i,j}\}_{1\leq\;i,j\;\leq c}$ and $\{ g_{k,j}\}_{1\leq\;k\;\leq d,\;1\leq\;j\;\leq c}$ be generalized filter functions associated to the multiplicity function
$\mu,$ that satisfy appropriate conditions and let $\{\Phi_1,\Phi_2, \cdots, \Phi_c\}\;\subseteq\;\widehat{V_0}$ be generalized scaling functions constructed as described above.  Setting 
$$\Psi_k\;=\widehat{D}(\sum_{j=1}^c\;g_{k,j}\Phi_j),\;1\leq\;k\;\leq d,$$
the $\{\Psi_k\;|\;1\leq\;k\;\leq d\}$ form a NTF wavelet family for dilation by $N.$
\end{theorem}

Current research of the author, Baggett, Merrill, and P. Jorgensen involves a generalization of certain results of Baggett, Courter and Merrill.  Bratteli and Jorgensen have found a group (called the loop group) which acts freely and transitively on $Lip_1$ filter systems $\{m_0,m_1,\cdots,m_{N-1}\}$ for dilation by $N$ (such a collection of filter functions is called an {\bf $m$-system} by Bratteli and Jorgensen), and have shown that each $Lip_1\;m$-system can be used to construct a NTF wavelet family with $N-1$ elements. Just as $m$ systems have been shown to give rise to functions on $\T$ with values in  $N\times N$ unitary matrices, we have shown that each generalized filter system $(h_{i,j})$  and $(g_{k,j})$ coming from a GMRA of Baggett, Courter and Merrill gives rise to unitary-matrix valued functions on $\T,$ whose dimensions now vary with $x\in \T.$ We have also found an analogue of the  loop-group which acts on the generalized filter functions 
$(h_{i,j})$  and $(g_{k,j})$ to give new functions satisfying the same orthogonality relations  \cite{BJMP}.  We conjecture that operator algebraic methods in addition to the results of \cite{BJMP}to will allow us to 
generalize Theorem 5.13 of Baggett, Courter and Merrill stated above.

As we have seen, operator algebras and direct integrals have played key roles in abstract wavelet theory, and we hope they will help in our new research project as well.  


\bibliographystyle{amsalpha}

\end{document}